\documentstyle[amssymb,amstex,amscd,11pt]{article}
\newcommand{\PP }{\Bbb P}
\newcommand{\QQ }{\Bbb Q}
\newcommand{\CC }{\Bbb C}
\newcommand{\ZZ }{\Bbb Z}
\newcommand{\proof }{{\bf Proof \ \ }}
\newcommand{\codim }{{\text {\rm codim}}}
\newtheorem{theorem}{Theorem}
\newtheorem{proposition}{Proposition}

\newtheorem{corollary}{Corollary}

\begin{document}

\author{Bumsig Kim \\
UC-Berkeley}
\title{Quot schemes for flags and Gromov invariants for flag varieties}
\date{4-25-95\\
Revised 12-1-95}
\maketitle

\section{Introduction}

To define Gromov-Witten invariants arise in mirror symmetry, there are two
general rigorous methods so far \cite{RT}\cite{Ko}. In particular Kontsevich
introduced the notion of stable maps for a compactification of moduli spaces.
For Grassmannians, however, there is a natural compactification of the space $%
Mor_d({\PP}^1,Gr(n,r))$ of all holomorphic maps from ${\PP}^1$ to
Grassmannians with a given degree $d$ where $Gr(n,r)$ is the Grassmannian of
all rank $r$ quotient vector spaces of ${\CC}^n.$ We may see $Mor_d({\PP}%
^1,Gr(n,r))$ as the set of all rank $r,$ degree $d,$ quotient bundles of $%
{\CC}^n\otimes {\cal O}_{{\PP}^1}.$ It is not a compact space.
Hence we come to a
Grothendieck's Quot scheme $Quot_d(Gr(n,r)),$ the set of all rank $r,$
degree $d,$ quotient sheaves of ${\CC}^n\otimes {\cal O}_{{\PP}^1}$ \cite
{Gr}$.$ It is proven to be a smooth projective variety by S. A. Str\o mme
\cite{St}. Bertram and Franco-Reina used Grothendieck's Quot schemes for
Gromov-Witten invariants and quantum cohomology of Grassmannians
respectively \cite{Be}\cite{FR}. In this paper using Quot schemes and
a localization theorem we study Gromov-Witten invariants for partial flag
varieties. The strategy is the following. We extend A. Bertram's result of
Gromov-Witten invariants for special Schubert varieties of Grassmannians to
the case of partial flag varieties. To do so a Grothendieck's Quot scheme
is generalized for flags and proven to be an irreducible, rational, smooth,
projective variety following Str\o mme \cite{St}. The Hilbert schemes for
flags have already studied \cite{Se}. On a partial
flag manifold there is an action by a special linear group. It induces an
action on the Quot scheme for flags. There is another action on it by the
multiplicative group ${\CC}^{\times }.$ It is induced from the action $%
\CC ^{\times }$ on ${\PP}^1.$ The analogous action by ${\CC}^{\times }$
does not exist on the Kontsevich moduli space. These two actions are
commutative. These together give isolated fixed points. Using a localization
by action an explicit formula of Gromov invariant for special Schubert
classes with a certain condition is given.
Note that Kontsevich also uses torus actions on his moduli
space of stable maps \cite{Ko}.
In his case he has to deal with summation over trees
since the fixed subsets are rather complicated. For projective spaces the
formula derived in the sequel is shown to agree to the residue formula \cite
{Ki}. The author does not know how to directly
relate with the result of Givental and
Astashkevich-Sadov's computation \cite{GK}\cite{AS}. Now we state our
main results.

\bigskip

For given integers $s_0=0<s_1=n-r_1<\cdots <s_l=n-r_l<n=s_{l+1},$ a flag
variety $Fl:=F(s_1,s_2,...,s_l;n)$ is, by definition, the set of all flags
of complex subspaces $V_1\subseteq V_2\subseteq \cdots \subseteq V={\Bbb C}^n,
\;\dim V_i=s_i$. There are universal vector bundles  $S_i$ and
universal quotient bundles $Q_i$ over $Fl$
with fibers ${\Bbb C}^{s_i}$ and $\CC ^n/\CC ^{s_i}$ respectively.
We are interested in a moduli space, the set
$Mor_d({\Bbb P}^1,Fl)$ of all morphisms $%
\varphi $ from ${\Bbb P}^1$ to $Fl$ with $<{\Bbb P}^1,c_1(\varphi
^{*}Q_k)>=d_k$, $d=(d_1,d_2,...,d_l)$.
Since $Fl$ is the fine moduli space such that the associated flag functor
is equivalent to $Mor(\cdot ,Fl)$ where the image of the functor at a scheme
$S$ is the set of all flag quotient bundles of $V\otimes\cal{O}_S$ with
ranks $r_i$. From the point of view as above, $Mor_d(\PP ^1 ,Fl)$ is
realized as the set of
all flag quotient bundles $(F_1,...,F_l)$ of $V\otimes\cal{O}_{\PP ^1}$ with
rank $r_i$ and degrees $d_i$ and hence it can be compactifyed by collecting
flag quotient sheaves. More precisely,

\bigskip

\begin{theorem}\label{thm1}
There is a smooth compactification $fQuot_d(Fl)$
of $Mor_d({\Bbb P}^1,Fl)$.
The underlying set of $fQuot_d(Fl)$
is the set of all flag quotient sheaves $(F_1,...,F_l)$
of $V\bigotimes {\cal O}_{{\Bbb P}^1}$ over ${\Bbb P}^1$ where
rank of $F_i$ is $r_i$ and its degree is
$d_i.$ Over the irreducible, rational, projective variety
$fQuot_d(Fl)\times \PP ^1$ there are tautological bundles $\cal {E}_i$
and sheaves $\cal {Q}_i$, $i=1,...,l$. They form exact sequences
$$ 0\rightarrow \cal {E}_i\rightarrow V\otimes O_{\PP ^1\times fQuot_d(Fl)}
\rightarrow \cal {Q}_i\rightarrow 0.$$
The induced sheaf morphisms $\cal {E}_i\rightarrow \cal{Q}_{i+1} $
are identically zero.
\end{theorem}

This fine moduli space $fQuot_d(Fl)$ will give the Gromov-Witten invariants
defined in \cite{KM}\cite{Ko}.

\begin{theorem}\label{thm2}
Let $p_1,...,p_N$ be fixed N distinct points in
$\PP ^1$. For i=1,...,N, let $\alpha _i$
be integers in $\{s_1,...,s_{l+1}\}$ and
let $\beta _i$
be positive integers less than $s_{\alpha _i+1}-s_{\alpha _i-1}$.
Then the number of
morphism $\varphi $ from $\PP ^1$ to $Fl$
such that each $\varphi (p_i)$ is in each the Poinc\'{a}re
dual Schubert subvariety to the classes $c_{\beta _i}(S_{\alpha _i})$
and $<\PP ^1 ,c_1(\varphi ^*Q_k)>=d_k,$ $%
k=1,...,l,$ is well-defined and it is
\[ \int _{fQuot_d(Fl)}\wedge _i c_{\beta _i} (\cal {E}_{p_i}^i). \]
The integration is not depend on the choices of the
point $p_1,...,p_N$ in $\PP ^1$.
\end{theorem}

By the torus action on $fQuot_d(fl)$ induced from the  standard
$\CC ^{\times }$-action on $\PP ^1$ and the
standard $(\CC ^{\times})^n$ on $Fl$ one
can apply Bott's residue formula to the above integration to get
\begin{theorem}\label{thm3}
The integration in the theorem \ref{thm2} is
\[
\sum_{\text{all integers as in (\ref{integer})}}\frac{\prod _i
(\sigma _{\alpha _i}^{\beta _i})
\prod (\text{characters as in (\ref{tang1}))}}
{\prod (\text{characters as in (\ref{tang2}))}}.
\]
The notations will be explained as follows.
\end{theorem}

Consider a sequence of data by nonnegative integers $d_{i,j}$ and $a_{i,j}$:
\begin{equation}
(d_{1,1},a_{1,1};...;d_{1,s_1},a_{1,s_1})\cdots
(d_{l,1},a_{l,1};...;d_{l,s_l},a_{l,s_l})  \label{integer}
\end{equation}

such that $d_{i,j}-a_{i,j}\geq d_{i+1,j}-a_{i+1,j}\geq 0,\;a_{i,j}\geq
a_{i+1,j}$ and $\sum_{j=1}^{r_i}d_{i,j}=d_i.$

Set $b_{i,j}:=d_{i,j}-a_{i,j}.$ Then we consider

\begin{eqnarray}
(p-a_{i,j})\hbar +\lambda _{j^{\prime }}-\lambda _{j},\;\text{for }%
0\leq p\leq a_{i,j^{\prime }}-1,\;1\leq j,j^{\prime }\leq s_i,  \nonumber
\label{tang1} \\
(b_{i,j}-p)\hbar +\lambda _{j^{\prime }}-\lambda _{j},\;\text{for }%
0\leq p\leq b_{i,j^{\prime }}-1,\;1\leq j,j^{\prime }\leq s_i,  \label{tang1}
\\
(p-a_{i,j})\hbar +\lambda _{m}-\lambda _{j}.\;\text{for }0\leq p\leq
d_{i,j},\;1\leq j\leq s_i,\;s_i+1\leq m\leq n, \nonumber
\end{eqnarray}
and
\begin{eqnarray} \label{tang2} \\
  (p-a_{i,j})\hbar +\lambda _{j^{\prime }}-\lambda _{j},\;
\text{for }0
&\leq &p\leq a_{i+1,j^{\prime }}-1,\;1\leq j,\leq s_i,\;1\leq j^{\prime
}\leq s_{I+1},  \nonumber \\
(b_{i,j}-p)\hbar +\lambda _{j^{\prime }}-\lambda _{j},\;
 \text{for }0
&\leq &p\leq b_{i+1,j^{\prime }}-1,\;1\leq j,\leq s_i,\;1\leq j^{\prime
}\leq s_{i+1},  \nonumber \\
(p-a_{i,j})\hbar +\lambda _{m}-\lambda _{j},\;
 \text{for }0 &\leq &p\leq
d_{i,j},\;1\leq j\leq s_i,\;\;s_{i+1}+1\leq m\leq n.  \nonumber
\end{eqnarray}

Finally let
$\sigma ^k_i$ be the $k$-th elementary symmetry function of
$a_{i,j}\hbar +\lambda _{j}$, $j=1,...,s_i$.

\bigskip

{\bf Acknowledgments and remarks: }My special thanks goes to my advisor A.
Givental for his wonderful guide. I would also like to thank A. Bertram, R.
Hartshorne, M. Kontsevich, and S. A. Str\o mme for answering my questions.
After finishing a preliminary version of the paper
I learned there are more advanced
results by A. Bertram \cite{Be1} and I. Ciocan-Fontanine \cite{CF} in some
cases. But works seem complementary. I would like to express my thanks to
I. Ciocan-Fontanine for pointing out an error
in the preliminary version of the paper. By \cite{Be1}\cite{CF} the condition
$\beta _i < s_{\alpha _i+1}-s_{\alpha _i-1}$ could be omitted.

\section{flag-Quot schemes}
All schemes will be assumed to be algebraic schemes over
an algebraically closed field ${\bf k}$ of characteristic $0\;$and all
sheaves will be quasi-coherent.
The space of all rank $r,$ degree $d,$ quotients of trivial sheaf $V\otimes
{\cal O}_{{\Bbb P}^1}$---or equivalently, all subsheaves of rank $s=n-r$,
degree $-d$---is a smooth, rational, irreducible, projective variety $%
Qout(V\otimes {\cal O}_{{\Bbb P}^1}$, Hilbert polynomial $rx+r+d)$ \cite{Gr}%
\cite{St}. Let us denote by $Quot_d(s,V)$ the Quot scheme. It could be
considered as a compactification of the space $Mor_d({\Bbb P}^1,Gr(n,r))$ of
all degree $d$ holomorphic maps from the projective line ${\Bbb P}^1$ to the
Grassmannian $Gr(n,r)$ which is the set of all rank $r$ quotient spaces of $V$
\cite{FR}. It is equipped with a universal locally free sheaf ${\cal E}$ and
a universal quotient sheaf ${\cal Q}$ over ${\Bbb P}^1\times Quot$ with an
exact sequence $0\rightarrow {\cal E}\rightarrow V\bigotimes {\cal O}_{{\Bbb P%
}^1\times Quot}\rightarrow {\cal Q}\rightarrow 0.$ They are flat over $Quot.$
For the special Schubert varieties, Gromov invariants can be defined via
Quot schemes as enumerative invariants \cite{Be}\cite{FR}. To extend their
results to flag varieties (not necessary complete flag varieties), we will
construct Quot schemes for flags, a compactification of $Mor_d(\PP ^1, Fl)$.

\bigskip

Let $s=(s_1,...,s_l),\;s_{i+1}>s_i>0,$ and $d=(d_{1,}...,d_l)$ be
multi-indices of nonnegative integers. For this moduli problem,
first we introduce a moduli functor $F^s_d$.
A contravariant functor $F_d^s$ from the category of schemes to the
category of sets is defined to be: for a scheme $S$, $F_d^s(S):=$the set
of all flag subsheaves $(E_1,E_2,...,E_l,E%
_n=V\bigotimes {\cal O}_{{\Bbb P}^1\times S})$ over ${\Bbb P}^1\times S$ where
sheaf ${\cal E}_i$ is subsheaf of ${\cal E}_j$ if $i<j,$
sheaves are flat over $S$, and
the rank of ${\cal E}_i$ is $s_i$ and its degree is $-d\;$over ${\Bbb P%
}^1$.
The functor $F_d^s\,$can be defined by quotients in a more transparent way,
for different data $0\rightarrow E\rightarrow V\otimes\cal{O}_{\PP ^1\times S}$
could give the same data
$V\otimes \cal{O}_{\PP ^1\times S}\rightarrow F\rightarrow 0$
\cite{FR}.

Now we will show the functor is representable. First some notations;
for a sheaf $\cal F$
over $\PP ^1 \times S$, denote by $\cal {F}(m)$
$\cal{F}\otimes \pi ^*\cal {O}(m)$, where $\pi$
is the projection from $\PP ^1\times S$ to $\PP ^1$. For the
second projection from
$\PP ^1\times S$ to $S$ we will use $\pi _S$.
Taking advantage of the existence of Quot schemes,
the obvious candidate for the scheme representing the functor $F^s_d$ is
the appropriate subscheme of
$Quot_{d_1}(s_1,n)\times ...\times Quot_{d_l}(s_l,n)$.

Over the product of Quot schemes there are universal subsheaves $\cal {E}_i$
and quotient sheaves $\cal {Q}_i$ induced each from $Quot_{d_i}(s_i,n)$.
Define $fQuot_d(Fl) $ as the degeneracy loci of
\[
\bigoplus_{i=1}^{l-1}
\left( (\pi _{\Pi Quot_{d_i}(s_i,n)}){*}{\cal E}_i(m)\rightarrow V/%
(\pi _{\Pi Quot_{d_i}(s_i,n)})_{*}({\cal E}_{i+1}(m))\right)
\]
for any $m\geq \max_i\{d_i\}-1$. Note that
$(\pi _{\Pi Quot_{d_i}(s_i,n)})_{*}{\cal E}_i(m)$ and $%
V/(\pi _{\Pi Quot_{d_i}(s_i,n)})_{*}({\cal E}_{i+1}(m))$ are locally free
so that there is no problem giving a scheme structure on $fQuot_d(Fl)$.

For shorthand, we will write $Quot_{d_i}(s_i,n)$ by $Quot_i$ and $%
fQuot_d(Fl)$ by $fQuot$, then now we are ready for

\begin{proposition}
The functor $F_d^s$ is representable by the (unique) projective scheme, $%
fQuot_d(Fl)$.
\end{proposition}

\proof
The statement means that for any scheme $S$, $F_d^s(S)=Mor(S,fQuot)$
in a functorial way. Let $(E_1,...,E_l)\in F_d^s(S).$ Then we have the
morphism $g:S\rightarrow Quot_1\times ...\times Quot_l$ from
the fine moduli property of
the Quot schemes.
We see that
$g^{*}(\pi _{\Pi Quot_i})_{*}{\cal E}_i(m)
\cong (\pi _S)_{*}E_i(m)$ naturally \cite{Mu}. The
fact that $E_i\subset E_{i+1}$ implies that $(\pi _S)_{*}E_i(m)\subset (\pi
_S)_{*}E_{i+1}(m).$ Hence
$g^{*}(\pi _{\Pi Quot_i})_{*}{\cal E}_i(m)\subset g^{*}
(\pi _{\Pi Quot_i})_{*}%
{\cal E}_{i+1}(m)$ and $g$ factor through $fQuot$.
$\Box $

\bigskip

Over $fQuot$ there are exact sequences of universal sheaves
\[
0\rightarrow {\cal E}_i\rightarrow V_{{\Bbb P}^1\times fQuot}\rightarrow
{\cal Q}_i\rightarrow 0
\]
and surjections ${\cal Q}_i\rightarrow {\cal Q}_{i+1}.$ Each ${\cal E}_i$
are locally free since it is a subsheaf of a locally free sheaf over ${\Bbb P}%
^1\times fQuot$ and it is flat over $fQuot.$ For given $f\in Mor(S,fQuot),$
the corresponding quotient sheaves over ${\Bbb P}^1\times S$ is just the pull
back $(id\times f)^{*}({\cal Q}_i).$

\subsection{Irreducibility and Smoothness}

To show the Quot scheme for flags is irreducible and smooth, one can simply
adapt Str\o mme's proof \cite{St}.

\begin{theorem}\label{thm4}
$fQuot$ is an irreducible, rational, nonsingular, projective
variety.
\end{theorem}

\proof
We will work by the language of quotients rather than subsheaves.
For $m=0,-1$ and $i=1,...,l,$ let
${\cal Q}_m^i$ be $(\pi _{fQuot})_{*}{\cal Q}_i(m)$,
a locally free sheaf over $%
fQuot$ of rank $(m+1)r_i+d_i$, and let
$X_m^i\rightarrow fQuot$ be the associated
principal $GL((m+1)r_i+d_i)$-bundle. One has a smooth morphism $\rho :$%
\[
\begin{tabular}{l}
$\prod_{fQuot,i=1}^l(X_{-1}^i\times _{fQuot}X_0^i)=:Y$ \\
$\;\;\;\;\;\;\;\;\downarrow \rho $ \\
$\;\;\;\;\;\;\;fQuot.$%
\end{tabular}
\]

We will show that $Y$ is an irreducible and smooth variety after finding an
isomorphism to a smooth irreducible affine quasi-variety. Since $\rho $ is
smooth, we conclude $fQuot$ is a smooth, irreducible, projective variety. Let
$%
N_m^i:=V^{(m+1)r_i+d_i}$ for $m=0,-1$ and$\;i=1,...,l.$ Here $V^r$ is, by
definition, the $r$-dimensional vector space over the ground field $k.$ Let $%
W:=Hom(V,V^{r_1+d_1})\times Hom_{{\Bbb P}^1}(\pi ^{*}V^{d_1}(-1),\pi
^{*}V^{r_1+d_1})$

$\times Hom(V^{d_1},V^{d_2})\times Hom(V^{r_1},V^{r_2+d_2})\times Hom_{{\Bbb P%
}^1}(\pi ^{*}V^{d_2}(-1),\pi ^{*}V^{r_2+d_2})$

$\times \cdots \times Hom(V^{d_{l-1}},V^{d_l})\times
Hom(V^{r_l},V^{r_l+d_l})\times Hom_{{\Bbb P}^1}(\pi ^{*}V^{d_l}(-1),\pi
^{*}V^{r_l+d_l}).$ Let $\overline{X}:=$associated affine space. On ${\Bbb P}^1%
 \times \overline{X},$ there is a tautological diagram:
\begin{eqnarray*}
&&
\begin{array}{ccc}
\pi _{\overline{X}}^{*}V^{d_l}(-1) & {\rightarrow }_{v_1} & \pi _{%
\overline{X}}^{*}V^{r_l+d_l} \\
&  & \uparrow \mu _l
\end{array}
\\
&&\,\,\,\,\,\,\,\,\,\,\,\,\,\,\,\,\,\,\,\,\,\,\,\,\,\,\,\,\,\,\,\,\,\,\,\,\,%
\,\,\,\,\,\,\,\,\vdots \\
&&
\begin{array}{ccc}
\pi _{\bar{X}}^{*}V^{d_2}(-1) & {\rightarrow }_{v_2} & \pi _{\bar{X}%
}^{*}V^{r_2+d_2} \\
&  & \uparrow \mu _2
\end{array}
\\
&&
\begin{array}{ccc}
\pi _{\overline{X}}^{*}V^{d_1}(-1) & {\rightarrow }_{v_1} & \pi _{%
\overline{X}}^{*}V^{r_1+d_1} \\
&  & \uparrow \mu _1 \\
&  & \pi _{\overline{X}}^{*}V.
\end{array}
\end{eqnarray*}
\

Let $Z\subset \overline{X}$ be the nonsingular irreducible quasi-variety
defined by the conditions

i) $v_i$ is injective on each fiber over $X$,

ii) the induced map $\overline{\mu }_i:\pi _Z^{*}V\rightarrow $coker ($%
v_i)$ is surjective, and

iii) coker$(v_i)$ is flat over $Z$ with rank $r_i$ and degree $d_i$ on the
fibers of $\pi _Z.$

The above conditions are open conditions in algebraic geometry. Let $%
W^{\prime }:=Hom(V^{d_1},V^{r_2+d_2})\times \cdots \times
Hom(V^{d_{l-1}},V^{r_l+d_l}).$ Let $X$ be a closed subvariety of $Z$ $\,$ of
codimension $\sum_{i=1}^ld_i(r_{i+1}+d_{i+1})$ defined by the inverse image
of a morphism
\begin{equation}
W\rightarrow W^{\prime }  \label{ww}
\end{equation}
measuring commutativity of the diagram
\[
\begin{array}{ccc}
\pi ^{*}_ZV^{d_{i+1}}(-1) & \rightarrow _d & \pi ^{*}_ZV^{r_{i+1}+d_{i+1}} \\
\uparrow a &  & \uparrow b \\
\pi ^{*}_ZV^{d_i}(-1) & \rightarrow _c & \pi ^{*}_ZV^{r_i+d_i}
\end{array}
\]
by $d\circ a-b\circ c.$

$X$ is a product of hypersurfaces defined by irreducible quadratic
polynomials. It is smooth away where the right vertical arrow is zero map.
And so $X$ is smooth and irreducible. There is a natural morphism $%
g:X\rightarrow fQuot$ by the construction of $X.$

By the construction,
$g^{*}{\cal Q}_m^i=(N_m^i)\otimes \cal {O}_X=((m+1)r_i+d_i){\cal O}_X.$
Therefore we have a morphism $s$ in the diagram
\[
\begin{tabular}{lll}
$X\ $ \  & ${\rightarrow }_s$ & ${%
\prod_{fQuot,i}(X_{-1}^i}\times _{fQuot}X_0^i)=:Y$ \\
$g\searrow $ &  & $\swarrow \rho $ \\
& $fQuot$ &.
\end{tabular}
\]

We will show $X$ and $Y$ are isomorphic finding the inverse of $s$ which
complete the proof.
We are given
isomorphisms on $Y$
\[
\lambda _m:\rho ^{*}{\cal Q}_m^i\rightarrow (N_m^i)_Y.
\]

By a proposition (1,1) in \cite{St},
we get a diagram on ${\Bbb P}^1\times Y$ except $\Uparrow$
\[
\begin{array}{ccc}
0\rightarrow \pi _Y^{*}N_{-1}^i(-1)\rightarrow & \pi _Y^{*}N_0^i &
\rightarrow (1\times \rho )^{*}{\cal Q}_i\rightarrow 0 \\
& \Uparrow & \uparrow \\
0\rightarrow \pi _Y^{*}N_{-1}^{i-1}(-1)\rightarrow & \pi _Y^{*}N_0^{i-1} &
\rightarrow (1\times \rho )^{*}{\cal Q}_{i-1}\rightarrow 0.
\end{array}
\]

Note here that$\;\pi _Y^{*}(\pi _Y)_{*}(1\times \rho )^{*}{\cal Q}_{i-1}=\pi
_Y^{*}$ $\rho ^{*}{\cal Q}_0^{i-1}\cong \pi _Y^{*}(N_0^{i-1})_Y.$

Since a morphism from free sheaf is determined and can be defined by a
morphism between the space of global sections, there exists a unique lifting
as indicated by the vertical arrow $\Uparrow $. By the defining property of $%
X$, there is an induced morphism which is the inverse of $s$.

Since $X$ is an irreducible smooth affine quasi-variety and $g$ is smooth of
relative dimension $\sum_{i=1}^l(d_i^2+(r_i+d_i)^2),$ $R$ is smooth and
irreducible. It's dimensions is $d_1(n-r_2)+d_2(r_1-r_3)+\cdots
+d_{l-1}(r_{l-2}-r_l)+d_lr_{l-1}$

$+nr_1+r_1r_2+r_2r_3+\cdots +r_{l-1}r_l-r_1^2-r_2^2-\cdots -r_l^2.$ The
rationality will be from Bialynski-Birula's theorem after considering an
action \cite{St}. The action will be studied in the following section 3.
$\Box $

\subsection{Gromov-Witten invariants and flag-Quot schemes}

 From now on we will work over the complex number field ${\Bbb C}$ to
consider complex manifolds.
We shall recall the definition of Gromov-Witten invariants for homogeneous
projective variety $X$.
The variety is always smooth. Denote by $%
\overline{{\cal M}}_N(X,d)$ the moduli stack of stable maps of degree $d$
and genus $0$. The stack is represented by
a smooth compact oriented orbifold. The same notations will be taken
for the stack and the coarse moduli orbifold.
It has morphisms, contraction $\pi ^X$ and
evaluations $ev_i$ at the $i$-th marked point:
\[
\begin{array}{ccc}
\overline{{\cal M}}_N(X,d) & {\rightarrow }_{ev_i}X &  \\
\downarrow _{\pi ^X} &  &  \\
\overline{{\cal M}}_N &  &
\end{array}
\]
where $\overline{{\cal M}}_N$ is the coarse moduli space of stable $n$
marked points of genus zero. $\overline{{\cal M}}_N$ is a smooth compact
oriented manifold.

The (tree level) Gromov-Witten classes $I_{N,d}^X:H^{*}(X, \QQ)^{\otimes
N}\rightarrow H^{*}(\overline{{\cal M}}_N,\QQ )$ are defined as follows:
\[
I_{N,d}^X(a_1\otimes \cdots \otimes a_N):=(\pi ^X)_{!}(ev_1^{*}(a_1)\otimes
\cdots \otimes ev_N^{*}(a_N)).
\]

In the sequel
we are interested in $I_{N,d}^X(a_1\otimes \cdots \otimes a_N)[p]
\in \QQ$ where $[p]$ is the homology class defined by a point $p$ in
$\overline{{\cal M}}_N$. If we choose any $N$ ordered distinct points $p_i$ in
$\PP ^1$, we can naturally embed $Mor_d(\PP ^1, X)$ into $(\pi ^X)^{-1}(p)$
for any generic point $p$. The boundary $(\pi ^X)^{-1}(p)\smallsetminus
Mor_d(\PP ^1, X)$
does not matter much, namely
\begin{proposition}
Let $Y_i$ be Schubert subvarieties of $X.$ Then $%
\bigcap_{i=1}^Nev_i^{-1}(g_iY_i)=\bigcap_{i=1}^Nev_i^{-1}(g_iY_i)\cap Mor_d(%
{\Bbb P}^1,X)$ for generic $g_i.$
\end{proposition}

\proof
The proof follows from the following general setting.%
$\Box$

\begin{proposition}
$M$ be open subvariety of a variety $\bar{M}.$
Suppose a connected algebraic group $G$ acts
transitively on another variety $X$.
Given a morphism $f:\bar{M}\rightarrow X$
subvarieties $Y_i$ of pure dimension, for generic $g_i\in G,$
\[
\bigcap_if^{-1}(g_iY_i)=(\bigcap_if^{-1}(g_iY_i))\cap M
\]
provided dimensional condition $\sum_i\codim Y_i=\dim\bar{M}$.

\end{proposition}

\proof
Apply Kleiman's theorem in Fulton's Book \cite{Fu}.
For generic $g_i$, $(\bar{M}\backslash M)\cap
\bigcap_if^{-1}(g_iY_i)=\emptyset $ and for generic $g_i,$ $%
(\bigcap_if^{-1}(g_iY_i))\cap M$ is proper. Hence for generic $g_i,$ $%
(\bigcap_if^{-1}(g_iY_i))\cap M$ is proper and ($\bar{M}\backslash M)\cap
\bigcap_if^{-1}(g_iY_i)=\emptyset $.
$\Box$

Since $Mor_d(\PP ^1, X) $ is a nonsingular quasi-projective variety, we can
consider its Chow group $A_*(Mor_d(\PP ^1,X))$ with products.
Let us use the same notation $ev_i$ for the restriction of the evaluation
map to $Mor_d(\PP ^1 ,X)$. For $[ev^{-1}_1(Y_1)]\cdot ...\cdot [ev^{-1}_1(Y_1)]
\in A_0(Mor_d(\PP ^1 ,X)$, in $\ZZ $ is
\[ \int _{Mor_d(\PP ^1 ,X)} [ev^{-1}_1(Y_1)]\cdot ...\cdot [ev^{-1}_1(Y_1)]
\] after summing up the coefficients of cycles of
points in $Mor_d(\PP ^1 ,X)$, which is well-defined in these intersections.
It is equal to
$I_{N,d}^X(a_1\otimes \cdots \otimes a_N)$ for
the Poincare dual classes $a_i$ of $Y_i$ because
$(\pi ^X)^{-1}(p)$ is a projective variety
and has a resolution of singularities to avoid the intersection theory of
algebraic (smooth) stacks.

We would like to do a similar thing in $Quot$ schemes following Bertram
\cite{Be}.

\begin{proposition}
(c.f. Bertram) Suppose $Y\subset X\,$ is an irreducible subvariety of
codimension $c$ and suppose $Z\subset Mor_d({\Bbb P}^1,X)$ is an irreducible
subvariety. Then for any $p\in {\Bbb P}^1$ and a generic translate $g,$ the
intersection $Z\cap ev_p^{-1}(gY)$ is either empty or has codimension $c$
in $Z$ where $ev_p$ denoted the evaluation map at $p$.
\end{proposition}

\proof
 Apply Kleiman's theorem in Fulton's Book \cite{Fu}. $\Box$

\begin{corollary}
\label{indexco}Let $c_i=$co$\dim _XY_i$ in the setting of the above
definition. If $\sum_{i=1}^Nc_i>\dim (Mor_d(C,X)),$ then, for generic
elements $g_1,...,g_N$ , $\bigcap_{i=1}^Nev_{p_i}^{-1}(g_iY_i)=\emptyset .$
If $\sum_{i=1}^Nc_i=\dim (Mor_d(C,X)),$ then, for generic elements $%
g_1,...,g_N$ , $\bigcap_{i=1}^Nev_{p_i}^{-1}(g_iY_i)$ is isolated or empty.
\end{corollary}
The points $p_i$ in the above could not be distinct.

\bigskip

Let $({\cal E}_i)_p$ be the restriction of the sheaf ${\cal E}_i$ at $p$ in $%
{\Bbb P}^1.$ Consider a commutative diagram
\[
\begin{array}{ccc}
Mor({\Bbb P}^1{\bf ,}Fl) & \rightarrow & Mor({\Bbb P}^1,Gr(n,r_i)) \\
\downarrow ev_p &  & \downarrow ev_p \\
Fl & \rightarrow & Gr(n,r_i)
\end{array}
\]
where $ev_p$ is the evaluation map at $p\in {\Bbb P}^1.$ Let $W$ be the
subspace of $V^{*}$ used for defining $Z,$ i.e., the special Schubert
varieties associated to $W$. Let $V_d(p,Z)$ be the degenerate locus of the
sheaf homomorphism
$W\bigotimes {\cal O}_{Quot}\rightarrow ({\cal E}_i)_p^{*}.$
In this setting we have

\begin{proposition}
$V_d(p,Z)$ represents the ($s_i+1-\dim (W))$-th Chern class of
$({\cal E}_i)_p^*$ over $fQuot$.
\end{proposition}

\proof
When the flag variety $Fl$ is a Grassmannian, it is
proven by A. Bertram \cite{Be}. For the general case, just consider the
morphism $fQuot\rightarrow Quot_i$ from the embedding $fQuot\rightarrow \prod
Quot_i$ followed by the projection $\prod Quot_i\rightarrow Quot_i$. It is
smooth since both schemes are smooth and the induced homomorphism between
tangent spaces is surjective after looking at \ref{ww}
in the proof of the theorem \ref{thm4}. This implies the degeneracy locus has
the expected dimension and $[V_d(p,Z)]$ in the Chow ring of the smooth
projective variety $fQuot$ is the $(\codim V_d(p,Z)$)-th Chern class
of $({\cal E}_i)_p^*$ which complete the proof.
$\Box$

Using a Pl\"{u}cker embedding and a stratification $Quot_d({\Bbb P}^1,{\Bbb P}%
^n)=\coprod_{0\leq m\leq d}C_m\times Mor_{d-m}({\Bbb P}^1,{\Bbb P}^n)$ (by
locally closed schemes) where $C_m$ is the $m$-th symmetric product of ${\Bbb %
P}^1,$ one can extend Bertram's result for flag varieties.

\begin{proposition}
\label{Ber}Let $Z_i$ be a special Schubert variety with $c_i$ codimension $%
\leq s_{k_{i+1}}-s_{k_{i-1}}-1$ representing a Chern class of $S_{k_i}^{*}.$
Suppose $\sum_{i=1}^Nc_i\geq \dim (fQuot),$ then $%
\bigcap_{i=1}^Nev_{d,p_i}^{-1}(g_iZ_i)=\bigcap_{i=1}^NV_d(p_i,g_iZ_i)$ for
distinct points $p_i\in {\Bbb P}^1$.
\end{proposition}

\proof
We will use induction on the total degree $%
|d|=d_1+\cdots +d_l.$ When $|d|=0,$ it can be done by Kleiman's theorem. Let
$\widetilde{C}_m=\prod_{i=1}^lC_{m_i},$ where $m=(m_1,...,m_l)$ is a
multi-index$.$ Using the Pl\"{u}cker morphisms on Quot schemes, let us
consider the morphism
\begin{eqnarray*}
J &:&fQuot\hookrightarrow \prod_{i=1}^lQuot_{d_i}(s_i,n)\rightarrow
\prod_{i=1}^lQuot_{d_i}({\Bbb P}^1{\Bbb ,P}^{M_i}) \\
&=&Mor_d({\Bbb P}^1 ,\prod_{i=1}^l{\Bbb P}^{M_i})\cup
\bigsqcup_{|m|=1}^{|d|}\widetilde{C}_m\times Mor_{d-m}({\Bbb P}^1 ,%
\prod_{i=1}^l{\Bbb P}^{M_i}).
\end{eqnarray*}

We would like to show that
\[
\bigcap_{i=1}^NV_d(p_i,g_iZ_i)\cap J^{-1}\left( \bigsqcup_{|m|=1}^{|d|}%
\widetilde{C}_m\times Mor_{d-m}({\Bbb P}^1 ,\prod_{i=1}^l{\Bbb P}%
^{M_i})\right) =\emptyset .
\]
Then we are done. To do so one has to show, for each $m$, $|m|>0$,
\begin{equation}
\emptyset =\bigcap_{i=1}^NV_d(p_i,g_iZ_i)\cap J^{-1}\left( \widetilde{C}%
_m\times Mor_{d-m}({\Bbb P}^1 ,\prod_{i=1}^l{\Bbb P}^{M_i})\right).
\label{inter}
\end{equation}

For any subset $P$ of $\{p_1,...,p_N\}$ let
\begin{eqnarray*}
A_P &=&\{\text{quotient sheaves }Q=(Q_1,...,Q_l)\in \text{LHS of (\ref{inter})
} \\
|\text{ }\dim _{k(p_i)}(Q_{k_i})_{p_i}\otimes k(p_i) &<&r_{k_i}\text{ iff }%
p_i\in P\}
\end{eqnarray*}
\[
A_P\subset \bigcap_{p_i\notin P}J^{-1}(\widetilde{C}_m\times
ev_{d-m,p_i}^{-1}(g_iZ_i)).
\]

But for generic $g_i,\;\bigcap_{p_i\notin P}J^{-1}(\widetilde{C}_m\times
ev_{d-m,p_i}^{-1}(g_iZ_i))=\emptyset $ since $\bigcap_{p_i\notin
P}ev_{d-m,p_i}^{-1}(g_iZ_i)=\emptyset $ by dimension counting in $%
fQuot(d-m;s_1,...,s_l;n)$:
\begin{eqnarray*}
\sum_{p_i\notin P}\codim (ev_{d-m,p_i}^{-1}(g_iZ_i)) &\geq &\dim
(fQuot_d(Fl)) \\
&&-\sum_{p_i\in P}(s_{k_{i+1}}-s_{k_{i-1}}-1) \\
&&\dim (fQuot_d(Fl)) \\
&&-\sum_{p_i\in P}(s_{k_{i+1}}-s_{k_{i-1}}). \\
&\geq &\dim (fQuot_{d-m}(Fl)).
\end{eqnarray*}

Since $p_i$ are distinct, we conclude the last inequality above.
$\Box$

The proof of the theorem \ref{thm2} follows from what are done.

\section{A Formula by Localization}

\subsection{Equivariant action on ${\cal E}_i\rightarrow fQuot\times {\PP}%
^1$}

By the standard action of $SL(n)\times PGL(2)$ on $V\times \PP ^1$,
the group acts on the space of stalks of $V\otimes \cal {O}_{fQuot\times
\PP ^1}$ and hence on the subsheaves $\cal {E}_i$ and $fQuot$. The action
on the sheaves is equivariant. In particular the maximal complex torus action
of $T\times \CC ^{\times}$ will formulate integrations of wedges
products of Chern classes of $(\cal {E}_i)_p$ as certain finite
sums of characters using the localization theorem \cite{AB}.

For simplicity of notations let us do it for the Quot schemes $Quot$.
Consider the action by $T\times\CC ^{\times}$
on $Quot\times {\PP}^1$, then the action
has a lift on the total space of the vector bundle ${\cal E}$. Let ${\cal E}%
_p$ be the restriction of the sheaf ${\cal E}$ at $p$ in ${\PP}^1$.
The action has the lifting to vector bundles $\cal{E}_0$ and $\cal{E}_\infty$.
It means $\cal{E}_{0(\infty )}$ is an equivariant vector bundle and its
equivariant Chern classes can be considered.
For other points, say $p$, transitive $%
PSL(2)$-action on ${\PP}^1$ will show $\cal{E}_0$, $\cal{E}_\infty$, and
$\cal{E}_p$ are isotropic:
\[
\begin{array}{ccc}
{\cal E} &  & {\cal E} \\
\downarrow &  & \downarrow \\
Quot\times {\PP}^1 & \rightarrow _g & Quot\times {\PP}^1
\end{array}
\begin{array}{ccc}
{\cal E}_0 &  & {\cal E}_p \\
\downarrow &  & \downarrow \\
Quot & \rightarrow _g & Quot
\end{array}
\]
where $g\cdot 0=p.$ In particular the Chern classes of ${\cal E}_p\,$ are
independent to $p$ since the map induced by $g$ is homotopic to identity.
Let $\frac 12\hbar \;($resp. $\lambda _i)$ is (are) the ${\CC}^{\times
}\;( $resp. $T)\;$characteristic classes. Then,
\begin{eqnarray}
\int_{Quot}\phi (c_{i_1}({\cal E}_{p_1}),...,c_{i_m}({\cal E}_{p_m}))
&=&\int_{Quot}\phi (c_{i_1}({\cal E}_0),...,c_{i_m}({\cal E}_0))  \nonumber
\\
&=&[\text{push forward of }\phi \text{ of equivariant classes }  \nonumber \\
&&\text{of }(c_1,...,c_r)\text{ at }{\cal E}_0]_{\hbar =\lambda _i=0}
\nonumber \\
&=&[\text{localization into components $P$ of the fixed subset},
\nonumber  \label{local} \\
&&\text{ i.e.,}\;\sum_P\text {push forward}
\frac{i_{*}^P\phi }{E(v_P)}]_{\hbar =\lambda
_i=0}  \label{local},
\end{eqnarray}
where $i^P$ is the inclusion $P\subset Quot$ and
$E(v_P) $ is the equivariant Euler class of
the normal bundle of P in $Quot$.
The last expression is independent to $\hbar $ and $\lambda _i$, without
letting them zeros, if the quasi-homogeneous degree of $\phi $ given
by degrees of the Chern classes agrees the
dimension of $Quot$. The complete analog hold for $fQuot$.

It is easy to see that the fixed subset consists of finite points. Therefore
$E(v_P)$ in (\ref{local})
is the equivariant Euler class of the normal space over the point
$P\in Quot.$ It is
the product of the complex characters of the representation of
$T\times\CC ^{\times}$
in the irreducible complex one dimensional subspaces of the tangent space of $%
Quot$ at $p.$ In the following subsection,
we devote ourselves to spell out the all fixed points
and all characteristics of the representation
to finish the proof of the theorem \ref{thm3}.

\subsection{Computation}

Let us use the standard maximal torus $T\times \CC ^{\times}$ in the
picture of flag manifolds $Fl$ and the projective line ${\PP}^1.$
Fix a sequence of $(e_{k_1},e_{k_2},...,e_{k_l})$ where $\{e_i\}_{i=1}^n$ is
the standard basis of $V.$ Then, for data (\ref{integer})
in introduction, one may associate
a flag of subsheaves
\[
{\cal O}(-d_{i,j})\longrightarrow {\cal O}(-d_{i+1,j})
\]
by the global section
\[
x^{(a_{i,j}-a_{i+1,j})}y^{(b_{i,j}-b_{i+1,j})}.
\]

It is a fixed point by the action. For such a (\ref{integer}) and a
sequence, we can associate any fixed point in $fQuot$. We have found all
fixed points.

Note that the tangent space at $x$ of a scheme $X$ is the first order
infinitesimal deformation $Mor_x(D,X),$ the set of all morphisms sending the
closed point of Spectrum of the ring $D$ of dual numbers
to $x$. Therefore, at a subsheaf ${\cal S}$ over ${\PP}$ of $%
V\otimes\cal{O}_{{\PP}^1},$ the tangent space of Quot schemes is the set of
flat
families of quotient sheaves over the Spec$D$ whose fiber over the closed
point of Spec$D$ is ${\cal S}$. It is
$Hom({\cal S},V\otimes\cal{O}_{{\PP}^1}/{\cal S}).$

For the flag-Quot scheme consider the following equivariant short exact
sequence at a fixed point ${\cal S}$ of $fQuot$
\begin{eqnarray*}
0 &\rightarrow &T_{{\cal S}}fQuot\left( d_1,...,d_l;s_1,...,s_l;n\right) \\
&\rightarrow &T_{{\cal S}}\{Quot_{d_1}(s_1,n)\times \cdots \times
Quot_{d_l}(s_l,n)\}\rightarrow \prod_{i=1}^lHom({\cal S}_i,{\cal Q}%
_{i+1})\rightarrow 0.
\end{eqnarray*}

At the fixed point associated to \ref{tang1} we find all characters of
irreducible subspace of $T_{{\cal S}}Quot$ by the torus action. They are,
for all $1\leq i\leq l,$
\begin{eqnarray*}
(p-a_{i,j})\hbar +\lambda _{j^{\prime }}-\lambda _{j},\;\text{for }0
&\leq &p\leq a_{i,j^{\prime }}-1,\;1\leq j,j^{\prime }\leq s_i, \\
(b_{i,j}-p)\hbar +\lambda _{j^{\prime }}-\lambda _{j},\;\text{for }0
&\leq &p\leq b_{i,j^{\prime }}-1,\;1\leq j,j^{\prime }\leq s_i, \\
(p-a_{i,j})\hbar +\lambda _{m}-\lambda _{j}.\;\text{for }0 &\leq &p\leq
d_{i,j},\;1\leq j\leq s_i,\;s_i+1\leq m\leq n,
\end{eqnarray*}

from
\begin{eqnarray*}
&&\bigoplus_{i=1}^lHom({\cal S}_i,{\cal Q}_i) \\
&=&\bigoplus_{i=1}^l[\bigoplus_{j=1,j^{\prime
}=1}^{s_i,s_i}Hom(.x^{a_{i,j}}y^{b_{i,j}}{\cal O}_{j}(-d_{i,j}),{\cal O}%
_{j^{\prime }}/x^{a_{i,j^{\prime }}}y^{b_{i,j^{\prime }}}\cal{O}_j) \\
&&\bigoplus \bigoplus_{j=1,m=s_i+1}^{s_j,n}Hom(x^{a_{i,j}}y^{b_{i,j}}{\cal O%
}_{j}(-d_{i,j}),{\cal O}_{m})].
\end{eqnarray*}

Characters from $\prod_{i=1}^{l-1}Hom({\cal S}_i,{\cal Q}_{i+1})$ are, for $%
1\leq i\leq l-1,$%
\begin{eqnarray*}
(p-a_{i,j})\hbar +\lambda _{j^{\prime }}-\lambda _{j},\;\text{for }0
&\leq &p\leq a_{i+1,j^{\prime }}-1,\;1\leq j,\leq s_i,\;1\leq j^{\prime
}\leq s_{i+1}, \\
(b_{i,j}-p)\hbar +\lambda _{j^{\prime }}-\lambda _{j},\;\text{for }0
&\leq &p\leq b_{i+1,j^{\prime }}-1,\;1\leq j,\leq s_i,\;1\leq j^{\prime
}\leq s_{i+1}, \\
(p-a_{i,j})\hbar +\lambda _{m}-\lambda _{j},\;\text{for }0 &\leq &p\leq
d_{i,j},\;1\leq j\leq s_i,\;\;s_{i+1}+1\leq m\leq n.
\end{eqnarray*}

The fiber space of ${\cal S}_i$ at the point $0$ has characters
\[
a_{i,j}\hbar +\lambda _{j}
\]
for $0\leq j\leq s_i.$ Therefore the $k$-th Chern character is the $k$-th
symmetric function in those characters. Let us denote it by $\sigma _i^k$.

The proof of theorem \ref{thm3} follows from the proposition 5.
\subsection{Projective spaces}

In this section we will relate the our result to the residue formula of
intersection pairing in \cite{Ki} for projective spaces. The
author does not know for the other cases.

Let $x$ be the Chern class of ${\cal O}_{{\PP}^n}(-1).\,$Then the
Gromov-Witten invariant $I^{\PP ^n}_{N,d}(x^{\otimes (n+1)d+n})$ is
\begin{equation}
\sum\Sb 0\leq i\leq n  \\ 0\leq k\leq d  \endSb \frac{(\lambda _i+k\hbar
)^{(n+1)d+n}}{\prod\Sb 0\leq p\leq d  \\ p\ne k  \endSb ((p-k)\hbar )\prod\Sb
0\leq q\leq d  \\ 0\leq j\ne i\leq n  \endSb ((q-k)\hbar +\lambda _j-\lambda
_i)}  \label{proj}
\end{equation}

\begin{proposition}
$\sum_{N=0}^{N=\infty }
\frac {1}{N!}q^d I^{\PP ^n}_{N.d}(x^{\otimes N})$ is a global residue
\[
\frac 1{2\pi }\oint \frac{f(x)dx}{x^{n+1}-q}
\] where $q$ is a formal variable.
\end{proposition}

\proof
The identity
\begin{eqnarray*}
\frac 1{2\pi }\oint \frac{x^{(n+1)d+n}dx}{x^{(n+1)(d+1)}} &=&\left[ \frac
1{2\pi }\oint \frac{x^{(n+1)d+n}}{\prod\Sb 0\leq i\leq n \\ 0\leq k\leq d
\endSb (x-\lambda _i-k\hbar )}\right] _{\lambda _i=\hbar =0} \\
&=&(\text{\ref{proj}})
\end{eqnarray*}
implies the proof.
$\Box$

e-mail address: {\it bumsig@@math.berkeley.edu}
\end{document}